\documentclass[12pt]{article}

\usepackage{amsfonts}
\usepackage{amssymb}
\usepackage{amsbsy}
\usepackage{url}
\usepackage[cmtip,arrow]{xy}
\usepackage{pb-diagram,pb-xy}

\def\CalC{{\mathcal{C}}}
\def\CalD{{\mathcal{D}}}
\def\CalS{{\mathcal{S}}}
\newtheorem{theorem}{Theorem}
\newtheorem{prop}{Proposition}
\newtheorem{cor}{Corollary}
\newtheorem{conj}{Conjecture}

\begin{document}
\title{\bf Symmetric Homology of Algebras}
\author{\bf S. Ault and Z. Fiedorowicz}
\date{August 11, 2007}
\maketitle

\begin{abstract}
\noindent
In this note, we outline the general development of a theory of symmetric homology
of algebras, an analog of cyclic homology where the cyclic groups are replaced by
symmetric groups. This theory is developed using the framework of crossed simplicial
groups and the homological algebra of module-valued functors. The symmetric homology
of group algebras is related to stable homotopy theory. Two spectral sequences for
computing symmetric homology are constructed. The relation to cyclic homology
is discussed and some conjectures and questions towards further work are discussed.

\end{abstract}

Symmetric homology is the analog of cyclic homology, where the cyclic groups
are replaced by symmetric groups. The second author and Loday \cite{FL} 
developed the notion of {\it crossed simplicial group \/} as a framework
for making this idea precise.

\bigskip

{\bf Definition 1} A {\it crossed simplicial group\/ } is a category $\Delta G$
whose objects are the sets $[n]=\{0, 1, 2,\dots,n\}$ for $n\ge 0$, which
contains the simplicial category $\Delta$, and such that any morphism
$[m]\to [n]$ factors uniquely as
$$[m]\stackrel{\cong}{\longrightarrow}[m]\stackrel{\gamma}{\longrightarrow}[n],$$
where $\gamma$ is a morphism in $\Delta$. The collection of groups
$\{G_n=Aut_{\Delta G}([n])^{op}\}_{n\ge0}$ are called the {\it underlying
groups\/} of $\Delta G$. The commutation relations implicit in $\Delta G$ endow
$\{G_n\}_{n\ge0}$ with the structure of a simplicial set (but not necessarily the
structure of a simplicial group).

\bigskip

The standard well-known example of a crossed simplicial group is $\Delta C$,
whose underlying groups are the cyclic groups $\{\mathbf{Z}_{n+1}\}_{n\ge0}$.
Less well-known is $\Delta S$, whose underlying groups are the symmetric
groups $\{\Sigma_{n+1}\}_{n\ge0}$. The description of $\Delta S$ given in \cite{FL} is difficult to work with.
A much nicer construction of this
category is due to Pirashvili \cite{Pir}.

\bigskip

{\bf Definition 2} The category $\Delta S$ has objects $[n]=\{0, 1, 2,\dots,n\}$.
A morphism $f:[m]\to[n]$ is a set function together with a specification of a
total order on the point preimages $\{f^{-1}(i)\}_{0\le i\le n}$. Composition
of morphisms $[m]\stackrel{f}{\longrightarrow}[n]\stackrel{g}{\longrightarrow}[p]$
is given by specifying the order on $(gf)^{-1}(i)=\amalg_{j\in g^{-1}(i)}f^{-1}(j)$,
as the block ordering specified by the ordering on $g^{-1}(i)$ and then
within each block by the ordering specified on each $f^{-1}(j)$. Any morphism
$f:[m]\to [n]$ decomposes uniquely as the permutation on $[m]$ specified by
$\amalg_{0\le i\le n}f^{-1}(i)$ followed by an order preserving function
$[m]\to[n]$, which is thus in $\Delta$. The cyclic category $\Delta C$ is the
evident subcategory of $\Delta S$.

\bigskip

{\bf Remark} Pirashvili's construction is a special case of a more general construction
due to May and Thomason \cite{MT}. This construction associates to any topological operad
$\{\CalC(k)\}_{n\ge0}$ a topological category $\widehat{\CalC}$ together with a functor
from this category to finite sets such that the inverse image of any function
$f:[m]\to[n]$ is the space $\prod_{i=0}^n\CalC(|f^{-1}(i)|)$. Composition in $\widehat{\CalC}$
is defined using the composition of the operad. May and Thomason refer to  $\widehat{\CalC}$
as the {\it category of operators\/} associated to $\CalC$. They were interested in the
case of an $E_\infty$ operad, but their construction evidently works for any operad. The
category of operators associated to the discrete $A_\infty$ operad $\mathcal{A}ss$, which
parametrizes monoid structures,  is precisely Pirashvili's construction of $\Delta S$.

\bigskip

Now given any small category $\CalC$ and any commutative ring $k$, one can define
homological algebra of covariant and contravariant functors
\mbox{$F:{\CalC}\longrightarrow k\mbox{-modules}$}. The simplest way to describe
this is to consider the ring $k[{\CalC}]$, which is the free $k$-module
generated by all the morphisms in $\CalC$. Multiplication is defined on
this basis by composition if the morphisms are composable and
0 otherwise. A covariant functor \mbox{$F:{\CalC}\longrightarrow k\mbox{-modules}$}
is then exactly the same thing as a left \mbox{$k[{\CalC}]$-module} structure
on $\bigoplus_{C\in\mbox{Obj}({\CalC})}F(C)$. Similarly contravariant functors
correspond to right $k[{\CalC}]$-modules (equivalently, left 
\mbox{$k[{\CalC}^\mathrm{op}]$-modules}).
One then defines for contravariant $F$ and covariant $G$,
$$\mbox{Tor}^{\CalC}_*(F,G) = 
\mbox{Tor}^{k[{\CalC}]}_*\left(\bigoplus_{C\in\mbox{Obj}({\CalC})}F(C),
\bigoplus_{C\in\mbox{Obj}({\CalC})}G(C)\right)$$
(There are some small technicalities that need to be checked, as the ring
$k[{\CalC}]$ does not have a multiplicative unit if $\CalC$ has infinitely
many objects. But it does have local units which are sufficient to carry
this out.)

If $A$ is a $k$-algebra, then the cyclic bar construction defines a functor
$B^{cyc}A: \Delta C^{op}\longrightarrow k\mbox{-modules}$, and cyclic
homology can be defined as
$$HC_*(A)=\mbox{Tor}^{\Delta C}_*(B^{cyc}A,\underline{k}),$$
where $\underline{k}:\Delta C\longrightarrow k\mbox{-modules}$ denotes the trivial functor
which takes every object to $k$ and every morphism to the identity.

However the results of \cite{FL} were discouraging as to the prospect of
an analogous definition of symmetric homology. First of all, the cyclic
bar construction does not extend to a functor 
$\Delta S^{op}\longrightarrow k\mbox{-modules}$. Secondly it was shown that
for any functor $F: \Delta S^{op}\longrightarrow k\mbox{-modules}$,
$\mbox{Tor}^{\Delta S}(F,\underline{k})$ is just the homology of the underlying simplicial
module of $F$, given by restricting $F$ to $\Delta^{op}$.

Subsequently, the second author \cite{F} noticed that the cyclic bar construction
extends not to a contravariant functor on $\Delta S$ but to a covariant
functor.

\bigskip

{\bf Definition 3} The {\it symmetric bar construction\/} is the functor
 $B^{sym}A: \Delta S\longrightarrow k\mbox{-modules}$ 
which takes the object $[n]$ to the $(n+1)$-fold tensor product $A^{\otimes n+1}$
of $A$ with itself over $k$. If $f:[m]\to[n]$ is a morphism in $\Delta S$,
then $B^{sym}(f)$ takes $a_0\otimes a_1\otimes a_2\otimes\dots\otimes a_m$ to
$b_0\otimes b_1\otimes b_2\otimes\dots\otimes b_n$, where
$b_i=\prod_{j\in f^{-1}(i)}a_j$, where the product is taken in the order
specified on $f^{-1}(i)$. 

\bigskip

The cyclic bar construction can be identified with
the composite
$$\Delta C^{op}\stackrel{D}{\cong}\Delta C\subset \Delta S
\stackrel{B^{sym}}{\longrightarrow}k\mbox{-modules},$$
where $D$ is a suitable duality isomorphism.

This now allows us to define symmetric homology as

\bigskip

{\bf Definition 4} ${\displaystyle HS_*(A) = \mbox{Tor}^{\Delta S}_*(\underline{k},B^{sym}A),}$
where $\underline{k}:\Delta S^{op}\longrightarrow k\mbox{-modules}$ denotes the trivial functor
which takes every object to $k$ and every morphism to the identity.

\bigskip

One can use the standard bar resolution of $\underline{k}$ to compute $HS_*(A)$ as the homology
of the simplicial abelian group $L_*(A)$, where
$$L_p(A) 
= \bigoplus 
k\left[[m_0]\stackrel{f_1}{\longrightarrow}
[m_1]\stackrel{f_2}{\longrightarrow}[m_2]\stackrel{f_3}{\longrightarrow}\dots
\stackrel{f_p}{\longrightarrow}[m_p]\right]\otimes A^{\otimes m_0+1}.$$
Here the direct sum ranges over all composable chains of morphisms in $\Delta S$ of length $p$.
The 0-th face consists of deleting $f_1$ and acting on $A^{\otimes m_0+1}$ via $B^{sym}(f_1)$.
The $p$-th face consists of dropping $f_p$. The other faces are given by composing $f_{i+1}$ with $f_i$.
The degeneracies consist of inserting identities.

If $A=k[M]$ is a monoid ring, then $HS_*(A)$ has a simpler description: it is the homology with $k$-coefficients
of the nerve of the category whose set of objects is the disjoint union $\amalg_{n\ge 0} M^{n+1}$. A
morphism from $(m_0,m_1,m_2,\dots,m_p)\in M^{p+1}$ to $(m'_0,m'_1,m'_2,\dots,m'_q)\in M^{q+1}$ is
a morphism $f:[p]\to[q]$ in $\Delta S$, such that $m'_i=\prod_{j\in f^{-1}(i)}m_j$. In the special
case when $M=J(X_+)$ is the free monoid on a generating set $X$ (for which we have $A=k[M]=T(X)$, the
tensor algebra on $X$) we have the following result. 

\begin{theorem} $HS_*(T(X)) = HS_*(k[J(X_+)]) = H_*(C_\infty(X_+);k)$, where $C_\infty$ denotes the
monad associated to the little $\infty$-cubes operad (\cite{May}, \cite{CLM}).
\end{theorem}

We may replace $C_\infty$ above by the monad associated to any $E_\infty$ operad.
In particular it is preferable to use the monad associated to the operad $\CalD$
(see Theorem 3 below).

If the monoid is a group $G$, we have the following result.

\begin{theorem} $HS_*(k[G]) = H_*(\Omega\Omega^\infty S^\infty(BG); k)$
\end{theorem}

The special case when $G$ is free abelian of rank $n$ is of particular interest. In this case
the group ring is the ring of Laurent polynomials in $n$ indeterminates. On the other hand
$BG$ is the $n$-torus which stably splits into a wedge of spheres. Thus we obtain

\begin{cor} $HS_*\left(k[t_1^{\pm},t_2^{\pm},\dots,t_n^{\pm}]\right)=
H_*(\Omega\Omega^\infty S^\infty(\bigvee_{i=1}^n\bigvee_{j=1}^{n!/(i!(n-i)!)}S^i);k)$
$${\textstyle \qquad\qquad=H_*(\prod_{i=1}^n\prod_{j=1}^{n!/(i!(n-i)!)}\Omega^\infty 
S^\infty(S^{i-1}); k)}$$
\end{cor}

Since the symmetric homology of the group completion of a commutative monoid is the group completion
of the symmetric homology of the monoid, a natural conjecture would be

\begin{conj}$HS_*(k[t_1,t_2,\dots,t_n])=$
$${\textstyle \qquad\qquad\qquad\qquad H_*\big(\prod_{i=1}^n C_\infty(S^0)\times
\prod_{i=2}^n\prod_{j=1}^{n!/(i!(n-i)!)}\Omega^\infty S^\infty(S^{i-1}); k\big)}$$
\end{conj}

In the case $n=1$, this conjecture is a special case of Theorem 1.

The $E_\infty$ structure visible in the above examples is a general phenomenon present in $HS_*(A)$
for any algebra. In order to make this precise, we need to enlarge the category $\Delta S$ by
adding an initial object $[-1]$.  Call the resulting enlarged category $\Delta S_+$, and let
$L_*^+(A)$ be the resulting enlarged bar complex. Then $\Delta S_+$ is a strict symmetric monoidal
category (with the monoidal structure given by the coproduct) and we have

\begin{theorem} (a) $HS_*(A)=H_*(L_*^+(A))$\newline
(b) $L_*^+(A)$ is an $E_\infty$ chain complex with respect to the action of the $E_\infty$ operad
$\CalD$.\newline
(c) If $k=Z_p$, $p$ a prime, then $HS_*(A)$ is equipped with Dyer-Lashof homology operations.
\end{theorem}

The $E_\infty$ chain operad which acts on $L^+_*(A)$ is the chain operad associated to the
operad $\CalD=\{E\Sigma_n\}_{n\ge0}$ which acts on strict symmetric monoidal (a.k.a. permutative)
categories \cite{May2}. This operad, in its simplicial form, is known as the Barratt-Eccles operad.

The following related result is joint work with Tomas Barros.

\begin{theorem} (a) The chain complex $L^+_*(A)$ is equipped with a Smith filtration (\cite{Berger},
\cite{Smith}). The $n$-stage of this filtration is an $E_n$ chain complex.\newline
(b) If $A=k[G]$ is a group ring, then the homology of the $n$-stage of the Smith filtration on
$L^+_*(A)$ is isomorphic to $H_*(\Omega^n S^{n-1}(BG);k)$.
\end{theorem}

While the chain complex $L_*(A)$ fortuitously lends itself to computations of $HS_*(A)$ in the
special cases of the monoid rings of free monoids and group rings, it is much too unwieldy
for computations in general, as it is infinite dimensional in each degree.  As a first step
in obtaining a more tractable chain complex, we have

\begin{prop} If $A$ is equipped with an augmentation $A\to k$ and $I$ denotes the augmentation
ideal, then the inclusion $L^{epi}_*(A)\subset L_*(A)$ is a chain homotopy equivalence, where 
$$L^{epi}_p(A) = \bigoplus k\left[[m_0]\stackrel{f_1}{\twoheadrightarrow}
[m_1]\stackrel{f_2}{\twoheadrightarrow}[m_2]\stackrel{f_3}{\twoheadrightarrow}\dots
\stackrel{f_p}{\twoheadrightarrow}[m_p]\right] \otimes
 B^{sym}_{m_0}I,$$
for $p>0$, where the $f_i$ are required to be epimorphisms. Here
$$B^{sym}_m I = \left\{
\begin{array}{ll}
A &\mbox{if } m=0\\
I^{\otimes m+1} &\mbox{if } m>0\\
\end{array}
\right.$$
Thus $HS_*(A)=H_*(L^{epi}_*(A))$.
\end{prop}

The chain complex $L^{epi}_*(A)$ in turn can be filtered in a couple of ways, giving rise to
spectral sequences for computing $HS_*(A)$. The simplest such spectral sequence arises by
filtering $L^{epi}_*(A)$ by the number of jumps: the $n$-th filtration of $L^{epi}_*(A)$ consists
of chains where at most $n$ of the $[m_{i-1}]\stackrel{f_i}{\twoheadrightarrow}[m_i]$ are strict (i.e.
$m_{i-1}>m_i$). We obtain

\begin{theorem}If $A$ is equipped with an augmentation with augmentation ideal $I$, then there is
a first quadrant spectral sequence converging to $HS_*(A)$ with
$$E^1_{p,q} = \bigoplus_{m_0>m_1>m_2>\dots>m_p\ge0} H_q\left(\Sigma_{m_p+1}^{op};
B^{sym}_{m_0}I\otimes 
k\left[\prod_{i=1}^p\mbox{Epi}_\Delta([m_{i-1}],[m_i])\right]\right)$$
\end{theorem}

Here $\mbox{Epi}_\Delta([m],[n])$ denotes the set of epimorphisms in $\Delta$ between $[m]$ 
and $[n]$.
The group homology is defined with respect to the group right action of $\Sigma_{m_p+1}$ 
given by the
isomorphism
$$\begin{array}{l}
\lefteqn{
B^{sym}_{m_0}I\otimes
k\left[\prod_{i=1}^p\mbox{Epi}_\Delta([m_{i-1}],[m_i])\right]
}\\

\qquad \cong B^{sym}_{m_0}I\otimes_{k[\Sigma_{m_0+1}]}k\left[\mbox{Epi}_{\Delta S}([m_0],[m_1])\right]
\otimes_{k[\Sigma_{m_1+1}]}k\left[\mbox{Epi}_{\Delta S}([m_1],[m_2])\right]\\
\qquad\qquad\qquad\otimes_{k[\Sigma_{m_2+1}]}\dots
\otimes_{k[\Sigma_{m_{p-1}+1}]}k\left[\mbox{Epi}_{\Delta S}([m_{p-1}],[m_p])\right]
\end{array}$$.

Here, the $B^{sym}_{m_0}I$ component comes before the chain of morphisms because we are
viewing it as a \textit{right} \mbox{$k\Sigma_{m_0+1}$-module} rather than a
\textit{left} \mbox{$k\Sigma_{m_0+1}^\mathrm{op}$-module}.

The differential $E^1_{p,q}\to E^1_{p-1,q}$ is an alternating sum of faces. The 0-th face takes
$$a_0\otimes a_1\otimes\dots a_{m_0}\otimes\left\{[m_0]\stackrel{f_1}{\twoheadrightarrow}
[m_1]\stackrel{f_2}{\twoheadrightarrow}[m_2]\stackrel{f_3}{\twoheadrightarrow}\dots
\stackrel{f_p}{\twoheadrightarrow}[m_p]\right\}$$
to
$$B^{sym}(f_1)(a_0\otimes a_1\otimes\dots a_{m_0})\otimes\left\{
[m_1]\stackrel{f_2}{\twoheadrightarrow}[m_2]\stackrel{f_3}{\twoheadrightarrow}\dots
\stackrel{f_p}{\twoheadrightarrow}[m_p]\right\}.$$
The middle faces compose consecutive arrows. The last face is induced by
$f_p^*:\Sigma_{m_p+1}\rightarrow\Sigma_{m_{p-1}+1}$, which is part of the simplicial structure
on the underlying groups $\{\Sigma_{n+1}\}_{n\ge0}$ of $\Delta S$.

Now, since the differential of $L^{epi}_*(A)$ reduces the filtration degree by at most
one, it can be shown that the differentials $E_{p,q}^r \to E^r_{p-r, q+r-1}$ must be 
trivial for $r \geq 2$.  Hence, the
spectral sequence collapses at the $E^2$ term.

This spectral sequence is still not very computationally useful as the $E^1$-term is
infinitely generated in each degree. A better spectral sequence is
obtained by filtering $L^{epi}_*(A)$ as follows:
$$F_mL^{epi}_p(A) =
\bigoplus_{m_0\le m} B^{sym}_{m_0}I\otimes k\left[[m_0]\stackrel{f_1}{\twoheadrightarrow}
[m_1]\stackrel{f_2}{\twoheadrightarrow}[m_2]\stackrel{f_3}{\twoheadrightarrow}\dots
\stackrel{f_p}{\twoheadrightarrow}[m_p]\right],$$
We obtain the following result.

\begin{theorem}If $A$ is equipped with an augmentation whose augmentation ideal $I$ is a free
$k$-module with basis $X$, then there is a spectral sequence converging strongly
to $HS_*(A)$ with
$$E^1_{p,q} = \bigoplus_{\overline{u}\in X^{p+1}/\Sigma_{p+1}}
\widetilde{H}_{p+q}(EG_{\overline{u}}\ltimes_{G_{\overline{u}}}N\CalS_p/N\CalS'_p;k)$$
\end{theorem}

Here $G_{\overline{u}}$ is the isotropy subgroup of the orbit $\overline{u}\in X^{p+1}/\Sigma_{p+1}$.
$N\CalS_p$ is the nerve of the category $\CalS_p$, which is defined as follows.  Let
$\{z_0,z_1,z_2,\dots\}$ be a countable set of indeterminates. First we define a larger category
$\tilde{\CalS}_p$. The objects of $\tilde{\CalS}_p$ are
formal tensor products $Z_0\otimes Z_1\otimes\dots\otimes Z_r$ where each $Z_i$ is a formal
(nonempty) product of the indeterminates $\{z_0,z_1,\dots,z_p\}$ so that 
$Z_0Z_1\dots Z_r=z_{\sigma(0)}z_{\sigma(1)}\dots z_{\sigma(p)}$ for some $\sigma\in\Sigma_{p+1}$.
In other words each $z_i$, $i=0,1,2,\dots,p$ occurs once and only once as a factor in exactly 
one of
the products $Z_j$, $j=0,1,2,\dots,r$. There is precisely one morphism in $\tilde{\CalS}_p$
$Z_0\otimes Z_1\otimes\dots\otimes Z_r\longrightarrow Y_0\otimes Y_1\otimes\dots\otimes Y_s$ 
iff each
$Y_i$ is a product of some of the monomials $Z_j$'s.  We then take $\CalS_p$ to be a skeletal 
subcategory of
$\tilde{\CalS}_p$. $\CalS_p$ is a poset. The nerve $N\CalS_p$ is contractible, since $\CalS_p$ contains
the initial object $z_0\otimes z_1\otimes\dots\otimes z_p$. We then take $\CalS'_p$ to be the 
subposet
obtained from $\CalS_p$ by deleting the initial object. Thus the quotient $N\CalS_p/N\CalS'_p$ has the
same homotopy type as the suspension of $N\CalS'_p$. The symmetric group $\Sigma_{p+1}$ acts on 
$\CalS_p$
by permuting the generators $\{z_0,z_1,z_2,\dots,z_p\}$. This induces an action on 
$N\CalS_p/N\CalS'_p$.
The differential $E^1_{p,q}\longrightarrow E^1_{p-1,q}$ is induced by the 0-th face map in 
$N\CalS_p$.

Thus a fundamental problem in computing symmetric homology is to determine the homotopy type of the
spaces $N\CalS_p/N\CalS'_p$ and to analyze the actions of the symmetric groups on these spaces. If $k$
is a field of characterisitic 0, just knowing the rational homology of these spaces and the action of
the symmetric groups on the homology would suffice to determine the $E^1$-term of the spectral sequence
of Theorem 6. However the chain complex of the simplicial nerve of $N\CalS_p/N\CalS'_p$ is too bulky
to permit computations except for very small values of $p$.

One can apply a similar technique, as is used to derive Theorem 5, to the nerve of the nonskeletal category
$\tilde{\CalS}_p$ to obtain a much smaller chain complex $Sym^{(p)}_*$, which computes the homology
of $N\CalS_p/N\CalS'_p$. The group of $i$-chains $Sym^{(p)}_i$ is the free abelian group on the objects of
$\tilde{\CalS}_p$ having the form $Z_0\otimes Z_1\otimes Z_2\otimes\dots\otimes Z_{p-i}$, modded 
out
by the equivalence relation generated by
\begin{eqnarray*}
\lefteqn{Z_0\otimes Z_1\otimes\dots\otimes Z_{j-1}\otimes Z_j\otimes\dots\otimes Z_{p-i}}\\
&= &(-1)^{(|Z_{j-1}|+1)(|Z_j|+1)}Z_0\otimes Z_1\otimes\dots\otimes Z_j\otimes Z_{j-1}\otimes
\dots\otimes Z_{p-i}
\end{eqnarray*}
where $|Z|$ denotes the length of the product. The boundary map in $Sym^{(p)}_*$ is an alternating sum
of faces, where each face consists of splitting a product $Z_j$ into a tensor product 
$Z'_j\otimes Z''_j$
(so that $Z_j=Z'_jZ''_j$ and the faces are ordered according to the position of the new $\otimes$. For example
$$\partial(z_2z_0z_3\otimes z_1z_4)= z_2\otimes z_0z_3\otimes z_1z_4-z_2z_0\otimes z_3\otimes 
z_1z_4 + z_2z_0z_3\otimes z_1\otimes z_4$$
The action of $\Sigma_{p+1}$ on $Sym^{(p)}_*$ is induced by permutation of the generators
$\{z_0,z_1,z_2,\dots,z_p\}$.

The direct sum $\bigoplus_{p\ge 0}Sym^{(p)}_*$ forms a bigraded differential algebra, where
$Sym^{(p)}_i$ is assigned bigrading $(p+1,i)$.  The product
$$\boxtimes: Sym^{(p)}_i\otimes Sym^{(q)}_j\longrightarrow Sym^{(p+q+1)}_{i+j}$$
is given by $Y\boxtimes Z=Y\otimes Z'$, where $Z'$ is obained from $Z$ by replacing each generator
$z_r$ by $z_{r+p+1}$ for $r=0,1,2,\dots,q$. The product is related to the boundary map by the
the relation
$$\partial(Y\boxtimes Z) = \partial(Y)\boxtimes Z + (-1)^i Y\boxtimes\partial(Z),$$ 
when $Y$ has bigrade $(p+1,i)$. Thus there is an induced map in homology:
$$\boxtimes: H_i(Sym^{(p)}_*)\otimes H_j(Sym^{(q)}_*)\longrightarrow H_{i+j}(Sym^{(p+q+1)}_*)$$
The product $\boxtimes$, both on the chain level and the homology level,  is not strictly skew commutative,
but rather skew commutative in a twisted sense:
$$Y\boxtimes Z = (-1)^{ij} \sigma Z\boxtimes Y$$
where $\sigma$ is the permutation which sends $0,1,2,\dots q$ to $p+1,p+2,\dots,p+q+1$ and
$q+1,q+2,\dots,p+q+1$ to $0,1,2,\dots,p$ in an order preserving way.

It is easy to compute the top degree homology groups. Let
\begin{eqnarray*}
b_p &= &z_0z_1z_2\dots z_p +(-1)^p z_1z_2\dots z_pz_0 +(-1)^{2p} z_2z_3\dots z_pz_0z_1\\
&& +\dots +(-1)^{p^2}z_pz_0z_1z_2\dots z_{p-1}.
\end{eqnarray*}
Then $b_p$ is a cycle and thus a homology class. As a $\mathbf{Z}[\Sigma_{p+1}]$-module,
$H_p(Sym^{(p)}_*)$ is generated by $b_p$ and as a representation $H_p(Sym^{(p)}_*)$ is
either the sign representation on $\mathbf{Z}_{p+1}$ (if $p$ is odd) or the trivial
representation on $\mathbf{Z}_{p+1}$ (if $p$ is even), induced up to $\Sigma_{p+1}$.
Thus $H_p(Sym^{(p)}_*)$ is free abelian of rank $p!$.

We summarize our calculations so far below.

\begin{theorem} For $p=0,1,2,3,4,5, 6, 7$ $H_*(Sym^{(p)}_*)$ are free abelian and have the
following Poincar\'e polynomials:
$$p_0(t)=1,\ p_1(t)=t,\ p_2(t)=t+2t^2,\ p_3(t)=7t^2+6t^3,$$
$$p_4(t)=43t^3+24t^4,\ p_5(t)=t^3+272t^4+120t^5,$$
$$p_6(t)=36t^4+1847t^5+720t^6,$$
$$p_7(t)=829t^5+13710t^6+5040t^7 $$
\end{theorem}

Ideally we would like to describe generators and relations for $\bigoplus_{p\ge 0}H_*(Sym^{(p)}_*)$
with respect to the module structures over the group rings of the symmetric groups and the
$\boxtimes$ product. The calculations summarized above show that besides the generators $b_i$
constructed above, there are additional generators in $H_3(Sym^{(4)}_*)$, $H_4(Sym^{(5)}_*)$,
$H_5(Sym^{(6)}_*)$, and $H_6(Sym^{(7)}_*)$. For now we only have
very limited understanding of these additional generators or of the relations between
the generators. For instance we have the following relation in $H_2(Sym^{(3)}_*)$
$$b_1\boxtimes b_1 = \left(1+[0312]+[1230]\right)b_2\boxtimes b_0,$$
where $[abcd]$ stands for the permutation $0\mapsto a, 1\mapsto b, 2\mapsto c, 3\mapsto d$.
The calculations also establish that $N\CalS_p/N\CalS'_p$ has the homotopy type of a wedge
of spheres for $p\le 6$. 

In recent work \cite{VZ} Vre\'cica and 
\v{Z}ivaljevi\'c  have connected $Sym^{(p)}_*$ to a certain well-studied
class of geometric complexes, known as chessboard complexes \cite{Wachs}. Using this they
have shown that
$Sym_*^{(p)}$ is $\lfloor \frac{2}{3}(p-1) \rfloor$-connected.
This result implies that the connectivity of the spaces $N\CalS_p/N\CalS'_p$
is an increasing function of $p$, hence the spectral sequence of Theorem 6
converges in the strong sense.  Indeed, for $m>\frac{3}{2}(i+1)$, there is an isomorphism
$$H_i(F_mL^{epi}_p(A))\stackrel{\cong}{\longrightarrow} HS_i(A)$$
If we denote by $\Delta^{(m)}S$ the full subcategory of $\Delta S$ consisting of the
objects $[0], [1], \ldots, [m]$, then it follows that for $m>\frac{3}{2}(i+1)$,
$$HS_i(A) = Tor^{\Delta^{(m)}S}_i(\underline{k},B^{sym}A).$$
Observe that $k\left[\Delta^{(m)}S\right]$ is a finite-dimensional unital ring, hence
if $A$ is finitely generated over a Noetherian ground ring $k$, then 
the increasing connectivity of the spaces
$N\CalS_p/N\CalS'_p$ implies that $HS_*(A)$ is finite dimensional over $k$ in each degree.
In the case when $A=k[G]$ is the group ring of a finite group, this also follows from Theorem 2, the
Atiyah-Hirzebruch spectral sequence for stable homotopy theory and Serre $\mathcal{C}$-theory.

Some questions are suggested by our partial computations of $H_*(Sym^{(p)}_*)$: Is it true that the
homology is always torsion-free? Or might it even be true that the spaces $N\CalS_p/N\CalS'_p$
are always wedges of spheres? Can the Vre\'cica and \v{Z}ivaljevi\'c connectivity result be
improved to $H_i(Sym^{(p)}_*)=0$ for $i\le p-r$, where
$$r=\left\lfloor\frac{\sqrt{8p+9}-1}{2}\right\rfloor\approx\sqrt{2p}?$$
If this were true, this would be the best possible connectivity result, since the sign
representation of $\Sigma_{p+1}$ has nontrivial multiplicity in all $H_i(Sym^{(p)}_*)$
for $p-r<i\le p$. The computed multiplicities of the trivial representations are also consistent
with this hypothesis.

\bigskip

We also have the following results on symmetric homology in degrees 0 and 1.

\begin{prop}(a) $HS_i(A)$ for $i=0,1$ is the homology of the following partial chain
complex
$$0\longleftarrow A \stackrel{\partial_1}{\longleftarrow} A\otimes A\otimes A
\stackrel{\partial_2}{\longleftarrow}(A\otimes A\otimes A\otimes A)\oplus A$$
where 
$$\partial_1(a\otimes b\otimes c) = abc -cba$$
$$\partial_2(a\otimes b\otimes c\otimes d) = ab\otimes c\otimes d + d\otimes ca\otimes b
+ bca\otimes 1\otimes d + d\otimes bc\otimes a,
\quad\partial_2(a) = 1\otimes a\otimes 1 
$$
(b) $HS_0(A) = A/[A,A]$ is the symmetrization of $A$ (as an algebra).
\end{prop}

We also have an elaboration of Theorem 1, which describes symmetric homology as the homology
of the $E_\infty$ symmetrization of an algebra. The idea is to simplicially resolve the algebra
by tensor algebras, then in each simplicial degree replace the tensor algebra by a free
$E_\infty$ chain algebra on the same generators, and finally to take the double complex
of the resulting simplicial chain algebra. A more precise formulation is

\begin{theorem} $HS_*(A)=H_*(B(D,T,A))$, where $B(D,T,A)$ is the 2-sided bar construction,
$T$ is the functor which takes a $k$-module to the tensor algebra on that module, $D$ is
the monad which takes a $k$-module to the free $\CalD$ chain algebra over that module (where
$\CalD$ is the same operad as in Theorem 3), and $B(D,T,A)$ is converted from a simplicial
chain complex to a double complex.
\end{theorem}

Finally we briefly discuss the relation between symmetric homology and cyclic homology. 
The relation between the cyclic bar construction and the symmetric bar construction, discussed
above, leads to a natural map 
$$HC_*(A)\longrightarrow HS_*(A).$$ 
The same analysis as in
Theorems 5 and 6 can be carried out for cyclic homology. The cyclic analog of $N\CalS'_p$
can be identified as a simplicial complex with the barycentric subdivision of  the
boundary of a $p$-simplex. The cyclic group acts on this $p-1$ sphere by cyclicly permuting
the vertices of the simplex. The cyclic analog of $N\CalS_p/N\CalS'_p$ is homotopy
equivalent to the suspension of this and is thus a $p$-sphere. One can then combine the
cyclic analog of Theorem 6 with the Serre spectral sequence for computing the homology of
the resulting half-smash products to obtain the standard spectral sequence for cyclic homology.

We can use the partial chain complex of Proposition 3 and an analogous one for cyclic homology
(c.f. \cite{Loday}, page 59) to describe the map $HC_i(A)\longrightarrow HS_i(A)$ for $i=0,1$.
These maps are induced by the following partial chain map:

\[
\begin{diagram}
  \node{ 0 }
  \node{ A }
  \arrow{w}
  \arrow{s,r}{\mathrm{id}}
  \node{ A \otimes A }
  \arrow{w,t}{ab - ba}
  \arrow{s,r}{a \otimes b \otimes 1}
  \node{A^{\otimes 3} \oplus A}
  \arrow{w,t}{\partial^C_2}
  \arrow{s,r}{f}\\
  \node{0}
  \node{A}
  \arrow{w}
  \node{A^{\otimes 3}}
  \arrow{w,t}{abc - cba}
  \node{A^{\otimes 4} \oplus A}
  \arrow{w,t}{\partial^S_2}
\end{diagram}
\]

The map $\partial^C_2$ takes $a \otimes b \otimes c \in A^{\otimes 3}$ to $ab \otimes c
- a \otimes bc + ca \otimes b$, and takes $a \in A$ to $1 \otimes a - a \otimes 1$.
The map $\partial^S_2$ is the map $\partial_2$ from \mbox{Proposition 3}.  $f$ is a map that is defined
on the first summand by
\begin{eqnarray*}
  a \otimes b \otimes c &\mapsto &a \otimes b \otimes c \otimes 1
  - 1 \otimes a \otimes bc \otimes 1 + 1 \otimes ca \otimes b \otimes 1\\
  &&+ 1 \otimes 1 \otimes abc \otimes 1 - b \otimes ca \otimes 1 \otimes 1
   - 2abc - cab
\end{eqnarray*}
and on the second summand by
$$a\mapsto 4a   - 1 \otimes 1 \otimes a \otimes 1$$

The map $HC_0(A)\longrightarrow HS_0(A)$ is the quotient map which takes the quotient of $A$ by
the $k$-module generated by all commutators onto the quotient of $A$ by the ideal
generated by all commutators. 

In a similar vein, Pirashvili and Richter (c.f. \cite{PR} and \cite{Sl}) have shown that
$$HC_*(A) = Tor^{\Delta S}_*(\underline{b},B^{sym}A),$$
where $\underline{b}$ is the contravariant functor on $\Delta S$ which is the cokernel
of $d_0-d_1: P_1\longrightarrow P_0$, where $P_i(-)=k[hom_{\Delta S}(-,[i])]$ and
$d_0-d_1$ induces the commutator map $a\otimes b\mapsto ab - ba$ on $B^{Sym}A$. Thus
the natural map $HC_*(A)\longrightarrow HS_*(A)$ is induced by the unique natural transformation
$\underline{b}\to\underline{k}$. Moreover the proof of Proposition 3 shows that $\underline{k}$
is the cokernel of $f-g:P_2\longrightarrow P_0$, where $f-g$ induces
$a\otimes b\otimes c\mapsto abc-cba$ on $B^{sym}A$, and the partial chain map above is induced
by a map from a partial projective resolution of $\underline{b}$ over $k[\Delta S]$ to a partial
projective resolution of $\underline{k}$ over $k[\Delta S]$.

The  Vre\'cica and \v{Z}ivaljevi\'c connectivity theorem, discussed above, implies that there is
a projective resolution of $\underline{k}$ which in degree $i$ is a finite direct sum of the projective
modules $P_m$ with $m\le\frac{3}{2}(i+1)$.

{DEPARTMENT OF MATHEMATICS, THE OHIO STATE UNIVERSITY, COLUMBUS, OH 43210-1174, USA\\
Email: ault@math.ohio-state.edu, fiedorow@math.ohio-state.edu}

\end{document}